\documentclass[a4paper,10pt]{amsart}

\usepackage[utf8]{inputenc} 
\usepackage{amsmath}
\usepackage{amsbsy}
\usepackage{amssymb}
\usepackage{amscd,amsthm}
\usepackage{graphicx}
\usepackage[mathcal]{eucal}
\usepackage{verbatim}
\usepackage[english]{babel}
\usepackage[dvigs]{epsfig}
\usepackage{dsfont}
\usepackage{longtable}
\usepackage{marvosym}
\usepackage{anysize}
\usepackage{hyperref}

\newcommand{\tl}{\text{li}}

\newcommand{\N}{\mathds{N}}
\newcommand{\R}{\mathds{R}}

\newcommand{\p}{\phantom}
\newcommand{\q}{\quad}

\newtheorem{thm}{Theorem}[section]
\newtheorem{lem}[thm]{Lemma}
\newtheorem{kor}[thm]{Corollary}

\newtheorem{prop}[thm]{Proposition}

\theoremstyle{definition}

\theoremstyle{remark}
\newtheorem*{rema}{Remark}

\title{On the sum of the first $n$ prime numbers} 
\author{Christian Axler}
\address{Institute of Mathematics\\ Heinrich-Heine-University Düsseldorf\\
40204 Düsseldorf, Germany}
\email{christian.axler@hhu.de}
\date{\today}

\begin{document}

\begin{abstract}
In this paper we establish a general asymptotic formula for the sum of the first $n$ prime numbers, which leads to a generalization of the most accurate 
asymptotic formula given by Massias and Robin. Further we prove a series of results concerning Mandl's inequality on the sum of the first $n$ prime numbers. We 
use these results to find new explicit estimates for the sum of the first $n$ prime numbers, which improve the currently best known estimates.
\end{abstract}

\maketitle

\section{Introduction}
Let $\pi(x)$ denote the number of primes not exceeding $x$. Hadamard \cite{hadamard1896} and de la Vall\'{e}e-Poussin \cite{vallee1896} independently proved a 
result concerning the asymptotic behavior for $\pi(x)$, namely $\pi(x) \sim \tl(x)$ as $x \to \infty$,
which is known as the \textit{Prime Number Theorem}.
In a later paper \cite{vallee1899}, where the existence of a zero-free region for the Riemann zeta-function $\zeta(s)$ to the left of the line $\text{Re}(s) = 
1$ was proved, de la Vall\'{e}e-Poussin also estimated the error term in the Prime Number Theorem by showing that
\begin{equation}
\pi(x) = \text{li}(x) + O(x e^{-a\sqrt{\log x}}), \tag{1.1} \label{1.1}
\end{equation} 
where $a$ is a positive absolute constant and the \emph{logarithmic integral} $\text{li}(x)$ is defined for every real $x \geq 0$ as
\begin{equation}
\text{li}(x) = \int_0^x \frac{dt}{\log t} = \lim_{\varepsilon \to 0+} \left \{ \int_{0}^{1-\varepsilon}{\frac{dt}{\log t}} + 
\int_{1+\varepsilon}^{x}{\frac{dt}{\log t}} \right \}. \tag{1.2} \label{1.2}
\end{equation}
Denoting the sum of the first prime numbers not exceeding $x$ by $S(x)$, Szalay \cite[Lemma 1]{szalay} used \eqref{1.1} to find
\begin{equation}
S(x) = \text{li}(x^2) + O ( x^2 e^{-a\sqrt{\log x}}). \tag{1.3} \label{1.3}
\end{equation}
Using \eqref{1.3} and integration by parts in \eqref{1.2}, we get the asymptotic expansion
\begin{equation}
S(x) = \frac{x^2}{2 \log x} + \frac{x^2}{4\log^2 x} + \frac{x^2}{4 \log^3 x} + \frac{3x^2}{8 \log^4 x} + O \left( \frac{x^2}{\log^5x} \right). \tag{1.4} 
\label{1.4}
\end{equation}
The first aim of this paper is to find explicit estimates for $S(x)$ in the direction of \eqref{1.4}. The current best such upper bound for $S(x)$ is due to 
Massias and Robin \cite[Th\'{e}or\`{e}me D(v)]{mr}. They found that $S(x) \leq x^2/(2\log x) + 3x^2/(10 \log^2 x)$ for every $x \geq 24\,281$. We start with 
the following result which improves the last inequality.

\begin{thm} \label{thm101}
For every $x \geq 110\,118\,925$, we have
\begin{displaymath}
S(x) < \frac{x^2}{2 \log x} + \frac{x^2}{4 \log^2 x} + \frac{x^2}{4\log^3 x} + \frac{5.3x^2}{8\log^4 x}.
\end{displaymath}
\end{thm}

The current best lower bound for $S(x)$ concerning \eqref{1.4} is also due to Massias and Robin \cite[Th\'{e}or\`{e}me D(ii)]{mr}. We find the following 
improvement.

\begin{thm} \label{thm102}
For every $x \geq 905\,238\,547$, we have
\begin{displaymath}
S(x) > \frac{x^2}{2 \log x} + \frac{x^2}{4 \log^2 x} + \frac{x^2}{4 \log^3 x} + \frac{1.2x^2}{8\log^4 x}.
\end{displaymath}
\end{thm}

Using an explicit estimate for $\tl(x^2)$, we find for the first time explicit bounds for the difference $S(x) - \tl(x^2)$ concerning \eqref{1.3} by 
establishing the following result.

\begin{thm} \label{thm103}
We have
\begin{displaymath}
- \frac{0.25x^2}{\log^4 x} < S(x) - \emph{li}(x^2) <  \frac{0.25x^2}{\log^4 x},
\end{displaymath}
where the left-hand side inequality is valid for every $x \geq 906\,484\,877$ and the right-hand side inequality holds for every $x \geq 110\,117\,797$.
\end{thm}

The case $x=p_n$, where $p_n$ denotes the $n$th prime number, is of particular interest. Here, $S(x) = \sum_{k \leq n}p_k$ is equal to the sum of the first $n$ 
prime numbers.
Massias and Robin \cite[p.\:217]{mr} found that
\begin{equation}
\sum_{k \leq n} p_k = \tl((\tl^{-1}(n))^2) + O(n^2 e^{-c\sqrt{\log n}}), \tag{1.5} \label{1.5}
\end{equation}
where $c$ is a positive absolute constant and $\tl^{-1}(x)$ is the inverse function of $\tl(x)$. Then they \cite[p.\:217]{mr} used \eqref{1.5} and a result of 
Robin \cite{robin1988} to derive the asymptotic expansion
\begin{equation}
\sum_{k \leq n} p_k = \frac{n^2}{2} \left( \log n + \sum_{i=0}^m \frac{A_{i+1}(\log \log n)}{\log^i n} \right) + O \left( \frac{n^2(\log \log 
n)^{m+1}}{\log^{m+1} n} \right), \tag{1.6} \label{1.6}
\end{equation}
where $m$ is a positive integer and the polynomials $A_k$ satisfy the formulas $A_0(x) = 1$ and $A_{k+1}' = A_k' - (k-1)A_k$. Unfortunately, this recursive 
formula for derivatives does not yield a description of the polynomials $A_k$, since the constant coefficient of the polynomials $A_k$ remains undetermined by 
this equation. We fix this problem by applying a method developed by Salvy \cite[Theorem 2]{salvy} and get the following theorem.

\begin{thm} \label{thm104}
Let $N$ be a positive integer. Then there exist uniquely determined monic polynomials $T_1, \ldots, T_{N-1}$ with real coefficients and $\emph{deg}(T_i) = i$, 
such that
\begin{displaymath}
\sum_{k \leq n} p_k = \frac{n^2}{2} \left( \log n + \log \log n - \frac{3}{2} + \sum_{i=1}^{N-1} \frac{(-1)^{i+1}T_i(\log \log n)}{i\log^i n} \right) + O \left( 
\frac{n^2(\log \log n)^N}{\log^Nn} \right).
\end{displaymath}
The polynomials $T_i$ can be computed explicitly. In particular,
\begin{itemize}
\item $T_1(x) = x - 5/2$,
\item $T_2(x) = x^2 - 7x + 29/2$,
\item $T_3(x) = x^3-12x^2+54x-185/2$,
\item $T_4(x) = x^4 - 52x^3/3 + 124x^2-442x + 1996/3$.
\end{itemize}
The polynomials $A_{i+1}$ given in \eqref{1.6} and polynomials $T_i$ are connected by the formula $T_i = (-1)^{i+1}iA_{i+1}$.
\end{thm}

The proof of Theorem \ref{thm104} is given in Section 5. The initial motivation for writing this paper was an inequality conjectured by Mandl concerning an 
upper bound for the sum of the first $n$ prime numbers, namely
\begin{equation}
\frac{np_n}{2} - \sum_{k \leq n} p_k \geq 0 \tag{1.7} \label{1.7}
\end{equation}
for every integer $n \geq 9$. This inequality originally appeared in \cite{rosser1975} without proof. In his thesis, Dusart \cite{pd} used the identity
\begin{displaymath}
\sum_{k \leq n} p_k = np_n - \int_2^{p_n} \pi(x) \, dx
\end{displaymath}
and explicit estimates for the prime counting function to prove that \eqref{1.7} indeed holds for every integer $n \geq 9$. The second goal of this paper is to 
study the sequence $(B_n)_{n \in \N}$, where $B_n$ denotes the left-hand side of \eqref{1.7}, in more detail.
For this purpose, we first derive an asymptotic expansion for $B_n$ by using a result of Cipolla \cite{cp} concerning an asymptotic expansion for the $n$th 
prime number. He proved that for every positive integer $N$ there exist uniquely determined monic polynomials $R_1, \ldots, R_{N-1}$ with real coefficients and 
$\deg(R_i) = i$, such that
\begin{displaymath}
p_n = n \left( \log n + \log \log n - 1 + \sum_{i=1}^{N-1} \frac{(-1)^{i+1}R_i(\log \log n)}{i\log^in} \right) + O\left( \frac{n (\log \log n)^N}{\log^Nn} 
\right).
\end{displaymath}
The polynomials $R_i$ can be computed explicitly. Setting $V_i = R_i - T_i$, where the polynomials $T_i$ are given by Theorem \ref{thm104}, we get the 
following asymptotic expansion for $B_n$.

\begin{thm} \label{thm105}
Let $N$ be a positive integer. Then,
\begin{displaymath}
B_n = \frac{n^2}{2}\left( \frac{1}{2} + \sum_{i=1}^{N-1} \frac{(-1)^{i+1}V_i(\log \log n)}{i\log^in} + O\left( \frac{(\log \log n)^N}{\log^Nn} \right) \right).
\end{displaymath}
The polynomials $V_s$ can be computed explicitly. In particular,
\begin{itemize}
\item $V_1(x) = 1/2$,
\item $V_2(x) = x-7/2$,
\item $V_3(x) = 3x^2/2 - 12x + 27$,
\item $V_4(x) = 2x^3 - 26x^2 + 124x - 221$.
\end{itemize}
\end{thm}


Since it is still difficult to compute $B_n$ for large $n$, we are interested in explicit estimates for $B_n$. From \eqref{1.7}, we get that $B_n > 0$ for every
integer $n \geq 9$. Hassani \cite[Corollary 1.5]{hassani2013} has found that the inequality $B_n > n^2/12$ holds for every integer $n \geq 10$. Up to now, the 
sharpest lower bound for $B_n$ is due to Sun \cite{sun}. He proved that $B_n > n^2/4$ for every integer $n \geq 417$. We improve Sun's result as follows.

\begin{thm} \label{thm106}
For every integer $n \geq 6\,309\,751$, we have
\begin{displaymath}
B_n > \frac{n^2}{4} + \frac{n^2}{4 \log n} - \frac{n^2(\log \log n - 2.9)}{4\log^2 n}.
\end{displaymath}
\end{thm}

In the other direction, we give the following explicit estimate for $B_n$, which improves the only known upper bound $B_n < 9n^2/4$, which holds for every 
integer $n \geq 2$, found by Hassani \cite[Corollary 1.5]{hassani2013}.

\begin{thm} \label{thm107}
For every integer $n \geq 256\,376$, we have
\begin{displaymath}
B_n < \frac{n^2}{4} + \frac{n^2}{4 \log n} - \frac{n^2(\log \log n-4.42)}{4\log^2 n}.
\end{displaymath}
\end{thm}

Theorem \ref{thm104} implies that
\begin{equation}
\sum_{k \leq n} p_k = \frac{n^2}{2} \left( \log n + \log \log n - \frac{3}{2} + \frac{\log \log n - 5/2}{\log n} - \frac{T_2(\log \log n)}{2\log^2 n} \right) +
O \left( \frac{n^2(\log \log n)^3}{\log^3n} \right), \tag{1.8} \label{1.8}
\end{equation}
where $T_2(x) = x^2 - 7x + 29/2$. We use the inequalities found in Theorems \ref{thm106} and \ref{thm107}, and combine them with some estimates for the 
$n$th prime number given in \cite[Theorems 1.1 and 1.4]{axler20173} to derive the following estimates for the sum of the first $n$ prime numbers, which 
refine the ones previously known. 

\begin{thm} \label{thm108}
For every integer $n \geq 1\,897\,700$, we have
\begin{displaymath}
\sum_{k \leq n} p_k < \frac{n^2}{2} \left( \log n + \log \log n - \frac{3}{2} + \frac{\log \log n - 5/2}{\log n} - \frac{(\log \log n)^2 - 7 \log \log n + 
13.567}{2\log^2 n} \right).
\end{displaymath}
\end{thm}

\begin{thm} \label{thm109}
For every integer $n \geq 2$, we have
\begin{displaymath}
\sum_{k \leq n} p_k > \frac{n^2}{2} \left( \log n + \log \log n - \frac{3}{2} + \frac{\log \log n - 5/2}{\log n}  - \frac{(\log \log n)^2 - 7 \log \log n + 
15.741}{2\log^2 n} \right).
\end{displaymath}
\end{thm}

\section{Proof of Theorem \ref{thm101}}

The following proof of Theorem \ref{thm101} is based on a recent obtained estimate for $\pi(x)$.

%
%
%

\begin{proof}[Proof of Theorem \ref{thm101}]
We denote the right-hand side of the required inequality by $f(x)$ and let $x_0 = 10\,166\,443\,802$. First, we consider the case $x \geq x_0$ and let $n = 
\pi(x) \geq 462\,277\,798$. We have
\begin{equation}
S(x) = \pi(p_n)p_n - np_n + \sum_{k \leq n} p_k. \tag{2.1} \label{2.1}
\end{equation}
Applying the upper bound for the prime counting function $\pi(x)$ given in \cite[Proposition 3.6]{axler20172} and the lower bound for $np_n - \sum_{k \leq n} 
p_k$ found in \cite[Theorem 1]{axler20174} to \eqref{2.1}, we get $S(x) < g(p_n)$, where
\begin{displaymath}
g(t) = \frac{t^2}{2 \log t} + \frac{t^2}{4 \log^2 t} + \frac{t^2}{4\log^3 t} + \frac{4.8t^2}{8\log^4 t} + \frac{4.5t^2}{4 \log^5 t} + \frac{28.5t^2}{8 \log^6 
t} + \frac{121.5t^2}{8\log^7 t} + \frac{25826.5t^2}{16\log^8 t}.
\end{displaymath}
Note that $g(t)$ is an increasing function for every $t \geq 17$. So we conclude that $S(x) < g(x)$ for every $x \geq x_0$. Since $g(t) < f(t)$ for every $t 
\geq x_0$, the proposition is proved for every $x \geq x_0$. A computer check shows that $f(p_i) \geq S(p_i)$ for every integer $i$ such that 
$\pi(110\,119\,007) \leq i \leq \pi(10\,166\,443\,802)$. Hence, $f(x) \geq S(x)$ for every $x$ with $110\,119\,007 \leq x \leq x_0$. Finally, we notice that 
$f(x) \geq S(x)$ for every $x$ satisfying $110\,118\,925 \leq x < 110\,119\,007$, which completes the proof.
\end{proof}

\begin{rema}
In \cite[Corollary 2.7]{delnic2}, Del\'{e}glise and Nicolas found a slightly weaker version of Theorem \ref{thm101}.
\end{rema}

\section{Proof of Theorem \ref{thm102}}

The currently best known lower bound for $S(x)$ is also due to Massias and Robin \cite[Th\'{e}or\`{e}me D(ii)]{mr}. They proved that the inequality $S(x) \geq
x^2/(2\log x) + 0.954x^2/(4 \log^2 x)$ holds for every $x \geq 70\,841$. In order to prove Theorem \ref{thm102}, we first note the following lemma, which can be 
found in \cite[Theorem 4.2]{ap}.

\begin{lem}[Abel's identity]  \label{lem301}
For any function $a : \mathds{N} \rightarrow \mathds{C}$ let $A(x) = \sum_{n \leq x} a(n)$, where $A(x)=0$ if $x<1$. Assume $g$ has a continuous derivative on
the interval $[y,x]$, where $0 < y < x$. Then we have
\begin{displaymath}
\sum_{y< n \leq x} a(n)g(n) = A(x)g(x) - A(y)g(y) - \int_{y}^{x}{A(t)g'(t)\ dt}.
\end{displaymath}
\end{lem}

The following proof of Theorem \ref{thm102} is based on the use of Lemma \ref{lem301} and some recently obtained estimates for the prime counting function.


\begin{proof}[Proof of Theorem \ref{thm102}]
First, we consider the case $x \geq 19\,027\,490\,297$. We denote the right-hand side of the required inequality by $f(x)$. Further, let $y=1, g(t) = t$ and
\begin{displaymath}
a(n) =
\begin{cases}
1 &\text {if $n$ is prime,} \\
0 &\text {otherwise.} \nonumber
\end{cases}
\end{displaymath}
We use Lemma \ref{lem301} to get
\begin{displaymath}
S(x) = \sum_{1<  n \leq x} a(n) g(n) = x\pi(x) - \int_1^x \pi(t) \, dt = x\pi(x) -  143 -  \int_{27}^x \pi(t) \, dt.
\end{displaymath}
Using the estimates for the prime counting function found in \cite[Propositions 3.6 and 3.12]{axler20172}, we see that
\begin{align*}
S(x) & > \frac{x^2}{\log x} + \frac{x^2}{\log^2 x} + \frac{2x^2}{\log^3 x} + \frac{5.85x^2}{\log^4 x} + \frac{23.85x^2}{\log^5 x} + \frac{119.25x^2}{\log^6 x} 
+ \frac{715.5x^2}{\log^7 x} + \frac{5008.5x^2}{\log^8 x} - 143 \\
& \p{\q\q} - \int_{27}^x \left( \frac{t}{\log t} + \frac{t}{\log^2 t} + \frac{2t}{\log^3 t} + \frac{6.15t}{\log^4 t} + \frac{24.15t}{\log^5 t} + 
\frac{120.75t}{\log^6 t} + \frac{724.5t}{\log^7 t} + \frac{6601t}{\log^8 t} \right) \, dt.
\end{align*}
Now, we apply \cite[Lemme 1.6]{pd} and \cite[Proposition 9]{axler2015} to this inequality, and get
\begin{align}
S(x) & > E_1 + \frac{26689x^2}{180\log x} + \frac{26689x^2}{360\log^2 x} + \frac{26689x^2}{360\log^3 x} + \frac{5327x^2}{48\log^4 x} + \frac{6661x^2}{30\log^5
x} + \frac{1663x^2}{3\log^6 x}  \tag{3.1} \label{3.1}  \\
& \p{\q\q} + \frac{3317x^2}{2\log^7 x} + \frac{10017x^2}{2\log^8x} - \frac{26599}{90}\,\tl(x^2), \nonumber
\end{align}
where $E_1$ is a constant with $E_1 \geq 111.708 > 0$. By \cite[Lemma 19]{axler2015}, we have
\begin{equation}
\tl(x^2) \leq \frac{x^2}{2\log x} + \frac{x^2}{4\log^2 x} + \frac{x^2}{4\log^3 x} + \frac{3x^2}{8\log^4 x} + \frac{3x^2}{4\log^5 x} + \frac{15x^2}{8\log^6 x} + 
\frac{45x^2}{8\log^7 x} + \frac{1575x^2}{64\log^8 x} \tag{3.2} \label{3.2}
\end{equation}
for every $x \geq 10^9$. Combined with \eqref{3.1}, we get
\begin{equation}
S(x) > \frac{x^2}{2 \log x} + \frac{x^2}{4 \log^2 x} + \frac{x^2}{4 \log^3 x} + \frac{3x^2}{20\log^4 x} + \frac{3x^2}{8\log^5 x} + \frac{3x^2}{16\log^6 x} - 
\frac{63x^2}{16\log^7 x} - \frac{289877x^2}{128\log^8x}, \tag{3.3} \label{3.3}
\end{equation}
which completes the proof for every $x \geq 19\,027\,490\,297$. To deal with the remaining case $905\,238\,547 \leq x < 19\,027\,490\,297$, we check with a 
computer that $S(p_i) \geq f(p_{i+1})$ for every integer $i$ with $\pi(905\,238\,547) \leq i \leq \pi(19\,027\,490\,297)$. Since $f'(x) > 0$ for every 
$x \geq 2.8$, we get $S(x) \geq f(x)$ for every $x \geq 905\,238\,547$.
\end{proof}

\begin{rema}
Recently, Theorem \ref{thm102} was independently found by Del\'{e}glise and Nicolas \cite[Corollary 2.7]{delnic2}.
\end{rema}

We obtain the following lower bound for $S(x)$, which corresponds to the first three terms of the asymptotic expansion \eqref{1.4}.

\begin{kor} \label{kor302}
For every $x \geq 152\,603\,617$, we have
\begin{equation}
S(x) > \frac{x^2}{2 \log x} + \frac{x^2}{4 \log^2 x} + \frac{x^2}{4 \log^3 x}. \tag{3.4} \label{3.4}
\end{equation}
\end{kor}

\begin{proof}
From Theorem \ref{thm102}, it follows that the required inequality holds for every $x \geq 905\,238\,547$. Similar to the proof of Theorem \ref{thm102}, we 
check \eqref{3.4} for smaller values of $x$ with a computer.
\end{proof}

The asymptotic formula \eqref{1.3} implies that
\begin{equation}
S(x) \geq \frac{x^2}{2 \log x} + \frac{x^2}{4 \log^2 x} \tag{3.5} \label{3.5}
\end{equation}
for all sufficiently large values of $x$. In 1988, Massias, Nicolas, and Robin \cite[Lemma 3(i)]{massias1989} proved that the inequality \eqref{3.5} holds for 
every $x$ such that $302\,791 \leq x \leq e^{90}$. Under the assumption that the Riemann hypothesis is true, Massias and Robin \cite[Th\'{e}or\`{e}me 
D(iv)]{mr} showed that the inequality \eqref{3.5} holds for every $x \geq 302\,971$. Further, they \cite[Th\'{e}or\`{e}me D(iv)]{mr} proved that the inequality 
\eqref{3.5} holds unconditionally for every $x$ such that $302\,971 \leq x \leq e^{98}$ and for every $x \geq e^{63864}$. Using Corollary \ref{kor302}, we 
fill this gap.

\begin{kor} \label{kor303}
The inequality \eqref{3.5} holds for every $x \geq 302\,971$.
\end{kor}

\begin{proof}
We only need to show that the desired inequality is valid for every $x$ such that $e^{98} < x < e^{63864}$. But this is a consequence of Corollary \ref{kor302}.
\end{proof}

%

\section{Proof of Theorem \ref{thm103}}

So far, we established explicit estimates for $S(x)$ in the direction of \eqref{1.4}. In the following proof of Theorem \ref{thm103}, where we establish for 
the first time explicit bounds for the difference $S(x) - \tl(x^2)$, we use an effective estimate for $\tl(x^2)$.


\begin{proof}[Proof of Theorem \ref{thm103}]
In the proof of Theorem \ref{thm101} it is shown that
\begin{displaymath}
S(x) < \frac{x^2}{2 \log x} + \frac{x^2}{4 \log^2 x} + \frac{x^2}{4\log^3 x} + \frac{4.8x^2}{8\log^4 x} + \frac{4.5x^2}{4 \log^5 x} + \frac{28.5x^2}{8 \log^6
x} + \frac{121.5x^2}{8\log^7 x} + \frac{25826.5x^2}{16\log^8x}
\end{displaymath}
for every $x \geq 10\,166\,443\,802$. By applying the corresponding lower bound for $\tl(x^2)$ given in \cite[Lemma 15]{axler2015}, we establish the 
correctness 
of the inequality
\begin{displaymath}
S(x) - \tl(x^2) <  \frac{0.225x^2}{\log^4 x} + \frac{0.375x^2}{\log^5x} + \frac{1.6875x^2}{\log^6 x} + \frac{9.5625x^2}{\log^7 x} +
\frac{1594.46875x^2}{\log^8 x}
\end{displaymath}
for every $x \geq 10\,166\,443\,802$. This completes the proof of the right-hand side inequality for every $x \geq 15\,884\,423\,625$. Similar to the 
proof of Theorem \ref{thm102}, we check with a computer that this inequality also holds for every $x$ such that $110\,117\,797 \leq x \leq 15\,884\,423\,625$. 
Analogously, we use \eqref{3.2}, \eqref{3.3} and a computer to verify that the desired left-hand side inequality is valid for every $x \geq 
906\,484\,877$.
\end{proof}

\begin{rema}
Under the assumption that the Riemann hypothesis is true, Del\'{e}glise and Nicolas \cite[Lemma 2.5]{delnic} improved \eqref{1.3} by showing that for every $x 
\geq 41$,
\begin{displaymath}
|S(x) - \tl(x^2)| \leq \frac{5}{24\pi}\, x^{3/2}\log x.
\end{displaymath}
\end{rema}


\section{Proof of Theorem \ref{thm104}}

In 1996, Massias and Robin \cite[p.\:217]{mr} found the currently the most accurate for the sum of the first $n$ primes, namely
\begin{equation}
\sum_{k \leq n} p_k = \frac{n^2}{2} \left( \log n + \log \log n - \frac{3}{2} + \frac{\log \log n - 5/2}{\log n} \right) + O \left( \frac{n^2(\log \log
n)^2}{\log^2 n} \right). \tag{5.1} \label{5.1}
\end{equation}
In Theorem \ref{thm104}, we give an asymptotic expansion for the sum of the first $n$ primes, which generalized \eqref{1.2}. The maintool for the given proof 
is a result of Salvy \cite[Theorem 2]{salvy}.


\begin{proof}[Proof of Theorem \ref{thm104}]
Let $N$ be a positive integer. We define
\begin{displaymath}
D_N(t) = \sum_{s=0}^N s!t^s.
\end{displaymath}
First, we note that repeated integration by parts in \eqref{1.2} gives
\begin{equation}
\tl(x) = \frac{x}{\log x} \left( D_N \left( \frac{1}{\log x} \right) + O \left( \frac{1}{\log^{N+1}x} \right) \right). \tag{5.2} \label{5.2}
\end{equation}
For $x > 1$, the logarithmic integral $\tl(x)$ is increasing with $\tl((1, \infty)) = \R$. Thus, we can define the inverse function $\tl^{-1} : \R \to (1,
\infty)$ by
\begin{equation}
\tl(\tl^{-1}(x)) = x. \tag{5.3} \label{5.3}
\end{equation}
The starting point of the proof is the asymptotic formula \eqref{1.5}. Using \eqref{5.2}, we get the asymptotic formula
\begin{equation}
\tl((\tl^{-1}(x))^2) = \frac{e^{2y}}{2y} \left( D_N \left( \frac{1}{2y} \right) + O \left(\frac{1}{y^{N+1}}  \right) \right) \tag{5.5} \label{5.5},
\end{equation}
where $y = \log \tl^{-1}(x)$. Next, we combine \eqref{5.2} and \eqref{5.3} to obtain $x = e^yy^{-1} D(1/y)$, where $D(t) = D_N(t) + O(t^{N+1})$. Now we apply
Theorem 2 of \cite{salvy} with $\alpha = 1, \beta = 2$, and $\gamma = -1$ to see that
\begin{displaymath}
\frac{e^{2y}}{2y} D_N \left( \frac{1}{2y} \right) = \frac{x^{2}\log x}{2} \sum_{i=0}^N \frac{Q_i( \log \log x)}{\log^ix} + O \left( \frac{x^2(\log \log
x)^N}{\log^Nx} \right),
\end{displaymath}
where the polynomials $Q_i \in \R[x]$ are defined by
\begin{equation}
Q_0 = 1, \; Q_{i+1}' = Q_i' - (i-1)Q_i. \tag{3.6} \label{3.6}
\end{equation}
Together with \eqref{1.5}, \eqref{5.5}, and the fact that $\tl^{-1}(x) \sim x \log x$ as $x \to \infty$, we conclude that
\begin{displaymath}
\sum_{k \leq n} p_k = \frac{n^{2}\log n}{2} \sum_{i=0}^N \frac{Q_i( \log \log n)}{\log^in} + O \left( \frac{n^2(\log \log n)^N}{\log^Nn} \right).
\end{displaymath}
By \eqref{3.6} and \eqref{5.1}, we have $Q_0(x) = 1$ and $Q_1(x) = x - 3/2$, respectively. Moreover, Theorem 2 of \cite{salvy} demonstrates how to compute the
value of the constant coefficient of the polynomials $Q_i$ for every integer $i$ with $2 \leq i \leq N$, which is not given by \eqref{3.6}. In the appendices 
of \cite{salvy}, one can find a Maple code for the computation of the polynomials $Q_2, \ldots, Q_N$ and it suffices to write
\begin{displaymath}
\texttt{(1/2)$\ast$theorem2}\_\texttt{part2(1,2,-1,D}\_\texttt{N(n),D}\_\texttt{N(n/2),n,N);}.
\end{displaymath}
Finally, we set $T_i = (-1)^{i+1}iQ_{i+1}$ for every integer $i$ satisfying $1 \leq i \leq N-1$. Then, \eqref{3.6} implies that the polynomials $T_i$ are monic 
with $\deg(T_i) = i$, which
completes the proof.
\end{proof}


\begin{rema}
The first part of Theorem \ref{thm104} was already proved by Sinha \cite[Theorem 2.3]{sin}.
\end{rema}

\section{Proof of Theorem \ref{thm105}}

Recall that $B_n = np_n/2 - \sum_{k \leq n} p_{k}$.
In this section, we use another result of Salvy \cite[Corollary 4]{salvy} (or Cipolla \cite{cp}) to give a proof of Theorem \ref{thm105} where we establish an 
asymptotic expansion for $B_n$.

\begin{proof}[Proof of Theorem \ref{thm105}]
Let $N$ be a positive integer. By Salvy \cite[Corollary 4]{salvy} (or Cipolla \cite{cp}) there exist uniquely determined monic polynomials $R_1, \ldots, R_{N-1}
$ with real coefficients and $\deg(R_i) = i$, so that
\begin{equation}
p_n = n \left( \log n + \log \log n - 1 + \sum_{i=1}^{N-1} \frac{(-1)^{i+1}R_i(\log \log n)}{i\log^in} \right) + O\left( \frac{n (\log \log n)^N}{\log^Nn}
\right). \tag{6.1} \label{6.1}
\end{equation}
Furthermore, in Appendix B.2 of \cite{salvy}, one can find a Maple code for the computation of the polynomials $R_1, \ldots, R_{N-1}$.
We set $V_i = R_i - T_i$ for every integer $i$ with $1 \leq i \leq N-1$, where the polynomials $T_i$ are given as in Theorem \ref{thm104}. Now it suffices to 
combine \eqref{6.1} and the asymptotic expansion given in Theorem \ref{thm104}.
\end{proof}



\section{Proof of Theorem \ref{thm106}}

In order to give a proof of Theorem \ref{thm106}, we first note the following proposition. Here, let
\begin{displaymath}
\gamma(n) = \frac{2.9\log^2n}{4\log^2 p_n} + \frac{\log^2 p_n \log^2 n + 16.7\log^2 n - \log^3 p_n \log n + \log^3p_n \log \log n}{4\log^3 p_n}.
\end{displaymath}

\begin{prop} \label{prop701}
For every integer $n \geq 6\,315\,433$, we have
\begin{displaymath}
B_n > \frac{n^2}{4} + \frac{n^2}{4 \log n} - \frac{n^2\log \log n}{4\log^2 n} + \frac{\gamma(n)n^2}{\log^2 n}.
\end{displaymath}
\end{prop}

\begin{proof}
First, we consider the case where $n \geq 440\,200\,309$. By \cite[Theorem 1]{axler20174}, we have
\begin{equation}
np_n - \sum_{k \leq n}p_k \geq \frac{p_n^2}{2 \log p_n} + \frac{3p_n^2}{4 \log^2 p_n} + \frac{7p_n^2}{4 \log^3 p_n} + L_1(n), \tag{7.1} \label{7.1}
\end{equation}
where $L_1(n) = (44.4p_n^2 \log^2 p_n + 184.2p_n^2 \log p_n + 937.5p_n^2)/(8\log^6 p_n)$. By \eqref{7.1} and the definition of $B_n$ it suffices to prove that
\begin{equation}
\frac{p_n^2}{2 \log p_n} + \frac{3p_n^2}{4 \log^2 p_n} + \frac{7p_n^2}{4 \log^3 p_n} + L_1(n) > \frac{np_n}{2} + \frac{n^2}{4} + \frac{n^2}{4 \log n} - 
\frac{n^2\log \log n}{4\log^2 n} + \frac{\gamma(n)n^2}{\log^2 n}. \tag{7.2} \label{7.2}
\end{equation}
For convenience, in the remaining part of the proof we write $p=p_n,y = \log n$, and $z = \log p$. It is easy to see that $937.5p^2 > 715.32npz + 
117.88n^2z^2$. We combine this inequality with the definition of $\gamma(n)$ to get
\begin{align*}
2n^2z^5y^2 & + 5.8n^2z^4y^2 + 55.5n^2z^2y^2(z-1.1) + 937.5p^2y^2 - 56.83n^2z^2y^2 - 715.32npzy^2\\
& > 2n^2z^6y - 2n^2z^6 \log y + 8\gamma(n)n^2z^6 + 22.1n^2z^3y^2.
\end{align*}
By Dusart \cite[Th\'{e}or\`{e}me 1.10]{pd}, we have $p > n(z - 1.1)$. Hence,
\begin{align*}
2n^2z^5y^2 & + 5.8n^2z^4y^2 + 184.2npzy^2(z-1.1) + 937.5p^2y^2 - 56.83n^2z^2y^2\\
& > 2n^2z^6y - 2n^2z^6 \log y + 8\gamma(n)n^2z^6 + 22.1n^2z^3y^2 + 128.7npz^2y^2 + 512.7npzy^2.
\end{align*}
Again, we use the inequality $p > n(z - 1.1)$ to obtain
\begin{align}
2n^2z^5y^2 & + 5.8n^2z^4y^2 + 184.2p^2zy^2 + 937.5p^2y^2 - 56.83n^2z^2y^2 \tag{7.3} \label{7.3} \\
& > 2n^2z^6y - 2n^2z^6 \log y + 8\gamma(n)n^2z^6 + 22.1n^2z^3y^2 + 128.7npz^2y^2 + 512.7npzy^2. \nonumber
\end{align}
Similar, we apply the inequality $p > n(z-1-1.15/z)$ found in \cite[Corollary 3.3]{axler20172} to \eqref{7.3} and see that
\begin{align*}
2n^2z^5y^2 + 8L_1(n)z^6y^2 & > 2n^2z^6y - 2n^2z^6 \log y + 8\gamma(n)n^2z^6 + 6n^2z^4y^2 + 10.3n^2z^3y^2 \\
& \p{\q\q} + 43.26n^2z^2y^2 + 32.6npz^3y^2 + 84.3npz^2y^2 + 461.64npzy^2.
\end{align*}
Analogously, we use the inequality $p > n(z-1-1/z-3.69/z^2)$ which is valid by \cite[Corollary 3.3]{axler20172} to get
\begin{align}
14p^2z^3y^2 + 8L_1(n)z^6y^2 & > 2n^2z^6y - 2n^2z^6 \log y + 8\gamma(n)n^2z^6 + 2n^2z^5y^2 + 2n^2z^4y^2 + 6.3n^2z^3y^2 \tag{7.4} \label{7.4} \\
& \p{\q\q} + 28.5n^2z^2y^2 + 10npz^4y^2 + 18.6npz^3y^2 + 70.3npz^2y^2 + 409.98npzy^2 \nonumber.
\end{align}
Next, we apply the inequality $p > n(z-1-1/z-3.15/z^2-14.25/z^3)$, see \cite[Corollary 3.3]{axler20172}, to \eqref{7.4} in a similar way to obtain
\begin{align*}
6p^2z^4y^2 + 14p^2z^3y^2 + 8L_1(n)z^6y^2 & > 2n^2z^6y^2 + 2n^2z^6y - 2n^2z^6 \log y + 8\gamma(n)n^2z^6 + 4npz^5y^2 \\
& \p{\q\q} + 4npz^4y^2 + 12.6npz^3y^2 + 51.4npz^2y^2 + 324.48npzy^2.
\end{align*}
Finally, by applying the inequality $p > n(z-1-1/z-3.15/z^2-12.85/z^3-81.12/z^4)$, which is fulfilled by \cite[Corollary 3.3]{axler20172}, we get
\begin{align*}
4p^2z^5y^2 & + 6p^2z^4y^2 + 14p^2z^3y^2 + 8L_1(n)z^6y^2 \\
& > 4npz^6y^2 + 2n^2z^6y^2 + 2n^2z^6y - 2n^2z^6 \log y + 8\gamma(n)n^2z^6.
\end{align*}
We divide the last inequality by $8 z^6y^2$ to obtain the inequality \eqref{7.2}, so the claim follows for every integer $n\geq 440\,200\,309$. We check 
the remaining cases with a computer.
\end{proof}

In 2012, Sun \cite{sun} proved that the inequality $B_n > n^2/4$ is valid for every integer $n \geq 417$. By proving Theorem \ref{thm106}, we improve Sun's 
lower bound for $B_n$.


\begin{proof}[Proof of Theorem \ref{thm106}]
For convenience, we write again $y = \log n$ and $z = \log p_n$. First, we consider the case where $n \geq 6\,315\,433$. By Proposition \ref{prop701} it 
suffices to show that $\gamma(n) \geq 2.9/4$. In \cite[p.\:42]{axler2013}, it is shown that for every $m \geq 255$,
\begin{equation}
\log m \geq 0.75 \log p_m. \tag{7.5} \label{7.5}
\end{equation}
Furthermore, we have $x^2 - 6.8x + 16.7 \cdot 0.75^2 > 0$ for every $x \geq 4.88$. Together with \eqref{7.5}, we get
\begin{equation}
16.7y^2 + (\log^2 y - (1 + 5.8)\log y)z^2 + 2.9z\log^2 y - 2.9z(\log y - 1) \geq 0. \tag{7.6} \label{7.6}
\end{equation}
From Dusart \cite[Proposition 5.15]{dusart2018} and the inequality $\log(1+t) \leq t$, which holds for every $t > -1$, follows that
\begin{equation}
z \leq y + \log y + \frac{\log y - 1}{y} + \frac{\log y - 2}{y^2}. \tag{7.7} \label{7.7}
\end{equation}
Using the result of Rosser \cite[Theorem 1]{rosser} that $p_m > m \log m$ for every positive integer $m$, we obtain
\begin{equation}
-z + \log y \leq - y.  \tag{7.8} \label{7.8}
\end{equation}
Hence, from \eqref{7.6}, we get
\begin{equation}
16.7y^2 + z^2( \log y - 1 - 2.9)\log y - 2.9zy\log y - 2.9z(\log y - 1) \geq 0. \tag{7.9} \label{7.9}
\end{equation}
Let $f(x) = 3.9(\log \log x - 2)/\log x$. Then it is easy to see that $f$ has a global maximum at $x_0 = 3$. Hence $f(x) \leq f(3) \leq 0.2$ for every $x > 1$.
Similary, we get $2.9(\log \log x - 1)/\log x \leq 0.4$ and $2.9(\log \log x - 2)/\log^2 x \leq 0.01$ for $x > 1$. Therefore,
\begin{displaymath}
\frac{3.9z^2(\log y - 2)}{y} + \frac{2.9z^2(\log y -1)}{y} + \frac{2.9z^2(\log y-2)}{y^2} < z^2.
\end{displaymath}
We combine this with \eqref{7.9} to obtain
\begin{align*}
z^2 & + 16.7y^2 + 2.9 zy^2 + z^2(\log y - 1 - 2.9)\log y\\
& \geq 2.9zy \left( y + \log y + \frac{\log y-1}{y} + \frac{\log y-2}{y^2} \right) + \frac{z^2(\log y-2)}{y} + \frac{2.9z^2}{y} \left(\log y-1 + \frac{\log
y-2}{y} \right).
\end{align*}
Now we use \eqref{7.7} to obtain
\begin{align*}
z^2 & + 16.7y^2 + 2.9 zy^2 + z^2(\log y - 1)\log y \\
& \geq 2.9z^2 \left( y + \log y + \frac{\log y-1 }{y} + \frac{\log y - 2}{y^2} \right) + \frac{z^2(\log y-2)}{y}.
\end{align*}
Again, by using \eqref{7.7}, we get
\begin{displaymath}
z^2 + 16.7y^2 + 2.9 zy^2 + z^2(\log y - 1)\log y \geq 2.9z^3 + \frac{z^2(\log y-2)}{y}
\end{displaymath}
and \eqref{7.8} implies
\begin{displaymath}
z^2y^2 + 16.7y^2 + 2.9 zy^2 - z^2y \left( y + \log y + \frac{\log y-1}{y} + \frac{\log y-2}{y^2} \right) + z^3 \log y \geq 2.9z^3.
\end{displaymath}
Finally we apply \eqref{7.7} to the last inequality and get $4z^3 \gamma(n) \geq 2.9z^3$.
Hence, the claim follows from Proposition \ref{prop701} for every $n \geq 6\,315\,433$. A computer check for smaller values of $n$ completes the proof.
\end{proof}


\section{Proof of Theorem \ref{thm107}}

We set
\begin{displaymath}
\kappa(n) = \frac{\log p_n \log^2 n + 4.1 \log^2 n - \log^2 p_n \log n + \log^2 p_n \log \log n}{4 \log^2 p_n} + \frac{r(\log p_n) \log^2 n}{8\log^6 p_n},
\end{displaymath}
where $r(x)$ is defined by
\begin{equation}
r(x) = 34.6x^3 + 207.1x^2 + 1431.56x + 28972.335, \tag{8.1} \label{8.1}
\end{equation}
to obtain the following proposition.

\begin{prop} \label{prop801}
For every integer $n \geq 256\,265$, we have
\begin{displaymath}
B_n < \frac{n^2}{4} + \frac{n^2}{4 \log n} - \frac{n^2\log \log n}{4\log^2 n} + \frac{\kappa(n)n^2}{\log^2 n}.
\end{displaymath}
\end{prop}

\begin{proof}
First, let $n \geq 841\,424\,976$; i.e. $p_n \geq 19\,033\,744\,403$. By \cite[Theorem 2]{axler20174} and the definition of $B_n$ it suffices to show that
\begin{equation}
\frac{np_n}{2} + \frac{n^2}{4} + \frac{n^2}{4\log n} - \frac{n^2\log \log n}{4\log^2 n} + \frac{\kappa(n)n^2}{\log^2n} > \frac{p_n^2}{2\log p_n} + 
\frac{3p_n^2}{4\log^2p_n} + \frac{7p_n^2}{4\log^3p_n} + U(n), \tag{8.2} \label{8.2}
\end{equation}
where
\begin{equation}
U(n) = \frac{45.6p_n^2}{8\log^4 p_n} +  \frac{93.9p_n^2}{4\log^5 p_n} + \frac{952.5p_n^2}{8\log^6 p_n} + \frac{5755.5p_n^2}{8\log^7 p_n} + 
\frac{116371p_n^2}{16\log^8 p_n}. \tag{8.3} \label{8.3}
\end{equation}
For convenience, we denote again $p=p_n,y = \log n$ and $z = \log p$. From the definiton of $\kappa(n)$ and $r(x)$, it follows that
\begin{align*}
2n^2z^8y & - 2n^2z^8\log y + 8\kappa(n)n^2z^8 + 3913.24n^2z^2y^2 \\
& = 2n^2z^7y^2 + 8.2n^2z^6y^2 + 34.6n^2z^5y^2 + 207.1n^2z^4y^2 + 1431.56n^2z^3y^2 + 32885.575n^2z^2y^2.
\end{align*}
By Rosser and Schoenfeld \cite[Corollary 1]{rosser1962}, we have $p < nz$. Hence, we obtain the inequality
\begin{align*}
2n^2z^8y & - 2n^2z^8\log y + 8\kappa(n)n^2z^8 + 3913.24n^2z^2y^2 + 25299.925npzy^2 \\
& > 2n^2z^7y^2 + 8.2n^2z^6y^2 + 34.6n^2z^5y^2 + 207.1n^2z^4y^2 + 1431.56n^2z^3y^2 + 58185.5npzy^2.
\end{align*}
Again, we use the inequality $p < nz$ to get
\begin{align*}
2n^2z^8y & - 2n^2z^8\log y + 8\kappa(n)n^2z^8 + 3913.24n^2z^2y^2 + 25299.925npzy^2 \\
& > 2n^2z^7y^2 + 8.2n^2z^6y^2 + 34.6n^2z^5y^2 + 207.1n^2z^4y^2 + 1431.56n^2z^3y^2 + 58185.5p^2y^2.
\end{align*}
Next, we apply the inequality $p < n(z-1)$, which was found by Dusart \cite{pd}, in a similar way to get
\begin{align}
2n^2z^8y & - 2n^2z^8\log y + 8\kappa(n)n^2z^8 + 610.47n^2z^3y^2 + 1871.21n^2z^2y^2 \tag{8.4} \label{8.4} \\
& \p{\q\q} + 3713.47npz^2y^2 + 19544.425npzy^2 \nonumber \\
& > 2n^2z^7y^2 + 8.2n^2z^6y^2 + 34.6n^2z^5y^2 + 207.1n^2z^4y^2 + 5755.5p^2zy^2 + 58185.5p^2y^2. \nonumber
\end{align}
The double usage of the inequality $p < n(z-1-1/z)$, see \cite[Corollary 3.9]{axler20172}, to \eqref{8.4} gives
\begin{align*}
2n^2z^8y & - 2n^2z^8\log y + 8\kappa(n)n^2z^8 + 110.4n^2z^4y^2 + 292.97n^2z^3y^2 + 1553.71n^2z^2y^2 + 635npz^3y^2  \\
& \p{\q\q} + 2760.97npz^2y^2 + 18591.925npzy^2 \\
& > 2n^2z^7y^2 + 8.2n^2z^6y^2 + 34.6n^2z^5y^2 + 952.5p^2z^2y^2 + 5755.5p^2zy^2 + 58185.5p^2y^2.
\end{align*}
Analogously, we apply the inequality $p < n(z-1-1/z-2.85/z^2)$, which was found in \cite[Corollary 3.9]{axler20172}, in a similar way to obtain
\begin{align*}
2n^2z^8y & - 2n^2z^8\log y + 8\kappa(n)n^2z^8 + 23.9n^2z^5y^2 + 51.9n^2z^4y^2 + 234.47n^2z^3y^2 + 1386.985n^2z^2y^2 \\
& \p{\q\q} + 129.3npz^4y^2 + 447.2npz^3y^2 + 2573.17npz^2y^2 + 18056.695npzy^2 \\
& > 2n^2z^7y^2 + 8.2n^2z^6y^2 + 187.8p^2z^3y^2 + 952.5p^2z^2y^2 + 5755.5p^2zy^2 + 58185.5p^2y^2.
\end{align*}
Next, we use that $p < n(z-1-1/z-2.85/z^2-13.15/z^3)$, see \cite[Corollary 3.9]{axler20172}, to get
\begin{align*}
2n^2z^8y & - 2n^2z^8\log y + 8\kappa(n)n^2z^8 + 6n^2z^6y^2 + 9.7n^2z^5y^2 + 37.7n^2z^4y^2 + 194n^2z^3y^2 + 1200.255n^2z^2y^2 \\
& \p{\q\q} + 31.4npz^5y^2 + 83.7npz^4y^2 + 401.6npz^3y^2 + 2443.21npz^2y^2 + 17457.055npzy^2 \\
& > 2n^2z^7y^2 + 8U(n) z^8y^2,
\end{align*}
where $U(n)$ is defined by \eqref{8.3}. Similar, we apply the inequality $p < n(z-1-1/z-2.85/z^2-13.15/z^3-70.7/z^4)$, which is valid by \cite[Corollary 
3.9]{axler20172}, to the last inequality and get
\begin{align*}
2n^2z^8y & - 2n^2z^8\log y + 8\kappa(n)n^2z^8 + 2n^2z^7y^2 + 2n^2z^6y^2 + 5.7n^2z^5y^2 + 26.3n^2z^4y^2 + 141.4n^2z^3y^2 \\
& \p{\q\q} + 917.455n^2z^2y^2 + 10npz^6y^2 + 17.4npz^5y^2 + 69.7npz^4y^2 + 361.7npz^3y^2 \\
& \p{\q\q} + 2259.11npz^2y^2 + 16467.255npzy^2 \\
& > 14p^2z^5y^2 + 8U(n) z^8y^2.
\end{align*}
Now, we use $p < n(z-1-1/z-2.85/z^2-13.15/z^3-70.7/z^4-458.7275/z^5)$, see \cite[Corollary 3.9]{axler20172}, in a analogical way to get
\begin{align}
2n^2z^8y^2 & + 2n^2z^8y - 2n^2z^8\log y + 8\kappa(n)n^2z^8 + 4npz^7y^2 + 4npz^6y^2 + 11.4npz^5y^2 \tag{8.5} \label{8.5} \\
& \p{\q\q} + 52.6npz^4y^2  + 282.8npz^3y^2 + 1834.91npz^2y^2 + 13714.89npzy^2 \nonumber\\
& > 6p^2z^6y^2 + 14p^2z^5y^2 + 8U(n)z^8y^2. \nonumber
\end{align}
Finally, we similary apply \cite[Theorem 3.8]{axler20172} to the inequality \eqref{8.5} to get
\begin{align*}
4npz^8y^2 & + 2n^2z^8y^2 + 2n^2z^8y - 2n^2z^8\log y + 8\kappa(n)n^2z^8 \\
& \p{\q\q} > 4p^2z^7y^2 + 6p^2z^6y^2 + 14p^2z^5y^2 + 8U(n)z^8y^2.
\end{align*}
We divide both sides of this inequality by $8z^8y^2$ to obtain the inequality \eqref{8.2} for every integer $n \geq 841\,424\,976$. We verify the 
remaining cases by using a computer.
\end{proof}


Now we use Proposition \ref{prop801} to give a proof of Theorem \ref{thm107}.

\begin{proof}[Proof of Theorem \ref{thm107}]
The proof consists of four steps. In the first step, we set $a_1 = 0.08$ and notice that
\begin{displaymath}
f(x) = 4a_1 ( x + \log x) + (x+4a_1 - \log x)\log \left( 1 + \frac{\log x - 1}{x} \right) - \log^2 x
\end{displaymath}
is positive for every $x \geq e^{19.63}$. In the following three steps, we write again $y = \log n$ and $z = \log p_n$, and consider the case $y \geq 19.63$. 
Then, $f(y) \geq 0$; i.e.,
\begin{equation}
\left( y + \log y + \log \left( 1 + \frac{\log y - 1}{y} \right) \right)(4a_1 + y - \log y) \geq y^2. \tag{8.6} \label{8.6}
\end{equation}
From Dusart \cite{dusart99} follows that
\begin{equation}
z \geq y + \log y + \log \left( 1 + \frac{\log y - 1}{y} \right). \tag{8.7} \label{8.7}
\end{equation}
We apply this inequality to \eqref{8.6} to get
\begin{equation}
8a_1z^8 \geq 2z^7y^2 - 2z^8y + 2z^8\log y. \tag{8.8} \label{8.8}
\end{equation}
In the third step, we set $a_2 = 1.025$ and $t(x) = 16a_2x^3\log x + 8a_2x^2\log^2 x - r(x)$, where $r(x)$ is defined by \eqref{8.1}. Then $t(x) \geq 0$ 
for 
every $x \geq 19.71$ and it follows that
\begin{displaymath}
16a_2z^5y^2\log z + 8a_2z^4y^2\log^2 z - r(z)z^2y^2 + (8a_2 - 8.2)z^6y^2 = z^2y^2t(z) \geq 0.
\end{displaymath}
The function $s \mapsto \log s/s$ is decreasing for every $s \geq e$. So, $\log(y)/y \geq \log(z)/z$ and we get
\begin{displaymath}
8a_2z^6(y + \log y)^2 - r(z)z^2y^2 - 8.2z^6y^2 \geq 0.
\end{displaymath}
By \eqref{8.7}, we obtain $z \geq y + \log y$. Hence $8a_2z^8 \geq r(z)z^2y^2 + 8.2z^6y^2$. Now, in the final step, we combine the last inequality with 
\eqref{8.8} to obtain
\begin{displaymath}
8.84 z^8 = 8(a_1+a_2)z^8 \geq 2z^7y^2 - 2z^8y + 2z^8\log y + r(z)z^2y^2 + 8.2z^6y^2 = 8\kappa(n)z^8.
\end{displaymath}
So, $\kappa(n) \leq 4.42/4$ for every integer $n \geq e^{19.63}$. We apply this to Proposition \ref{prop801}, which completes the proof for every integer $n 
\geq e^{19.63}$. We conclude by a direct computation.
\end{proof}

\begin{rema}
Theorem \ref{thm107} improves the only known upper bound $B_n < 9n^2/4$, which holds for every integer $n \geq 2$, found by Hassani \cite[Corollary 
1.5]{hassani2013}.
\end{rema}

\section{Proof of Theorems \ref{thm108} and \ref{thm109}}

In 1998, Dusart \cite{pd} proved that the inequality
\begin{displaymath}
\sum_{k \leq n} p_k \leq \frac{n^2}{2} \left( \log n + \log \log n - \frac{3}{2} + \frac{\log \log n - 2.29}{\log n} \right)
\end{displaymath}
holds for every integer $n \geq 10\,134$. In this section, we use the identity
\begin{equation}
\sum_{k \leq n} p_k = \frac{np_n}{2} - B_n, \tag{9.1} \label{9.1}
\end{equation}
the inequalities stated in Theorems \ref{thm106} and \ref{thm107}, and some explicit estimates for the $n$th prime number given in \cite{axler20173} to find 
proofs of Theorems \ref{thm108} and \ref{thm109}. We start with the proof of Theorem \ref{thm108}.

\begin{proof}[Proof of Theorem \ref{thm108}]
We combine \eqref{9.1}, \cite[Theorem 1.1]{axler20173}, and Theorem \ref{thm106}, to get that the required inequality holds for every integer $n \geq 
46\,254\,381$. The remaining cases are verified with a computer.
\end{proof}

Based on Theorem \ref{thm108}, we obtain the following upper bound for the sum of the first $n$ prime numbers, which corresponds to the first four terms of 
the asymptotic expansion found in Theorem \ref{thm104}.

\begin{kor} \label{kor603}
For every integer $n \geq 115\,149$, we have
\begin{equation}
\sum_{k \leq n} p_k < \frac{n^2}{2} \left( \log n + \log \log n - \frac{3}{2} + \frac{\log \log n - 5/2}{\log n} \right). \tag{9.2} \label{9.2}
\end{equation}
\end{kor}

\begin{proof}
Theorem \ref{thm108} implies the validity of \eqref{9.2} for every integer $n \geq 1\,897\,700$. It remains to check the required inequality for smaller 
values of $n$ with a computer.
\end{proof}

The current best lower bound for the sum of the first $n$ primes in the direction of \eqref{1.8} is also due to Dusart \cite[Lemme 1.7]{pd}. He proved that
\begin{displaymath}
\sum_{k \leq n} p_k \geq \frac{n^2}{2} \left( \log n + \log \log n - \frac{3}{2} \right)
\end{displaymath}
for every integer $n \geq 305\,494$. Using \cite[Theorem 1.2]{axler20173} and Theorem \ref{thm107}, we finally give the following proof of Theorem 
\ref{thm109}.


\begin{proof}[Proof of Theorem \ref{thm109}]
Applying \cite[Theorem 1.4]{axler20173} and Theorem \ref{thm107} to \eqref{9.1}, we get that the desired inequality holds for every integer $n \geq 
256\,376$. For the remaining cases, we use a computer.
\end{proof}

\section*{Acknowledgement}
I would like to thank Jean Pierre Massias and Guy Robin whose paper inspired me to deal with this subject. I would also like to thank R. for being a 
never-ending inspiration.


\end{document}